\documentclass[11pt,reqno]{amsart}
\usepackage{enumerate, latexsym, amsmath, amsfonts, amssymb, amsthm, color}
\def\pmod #1{\ ({\rm{mod}}\ #1)}
\def\Z{\Bbb Z}
\def\N{\Bbb N}

\def\Q{\Bbb Q}

\def\l{\left}
\def\r{\right}
\def\bg{\bigg}
\def\({\bg(}
\def\){\bg)}
\def\t{\text}
\def\f{\frac}

\def\mo{{\rm{mod}\ }}

\def\ord{{\rm ord}}

\def\ls{\leqslant}
\def\gs{\geqslant}

\def\sm{\setminus}

\def\ve{\varepsilon}

\def\eq{\equiv}

\def\da{\delta}

\def\Proof{\noindent{\it Proof}}

\theoremstyle{plain}
\newtheorem{theorem}{Theorem}

\newtheorem{lemma}{Lemma}
\newtheorem{corollary}{Corollary}
\newtheorem{conjecture}{Conjecture}
\theoremstyle{definition}

\theoremstyle{remark}
\newtheorem{remark}{Remark}

\allowdisplaybreaks

 \vspace{4mm}

\begin{document}

\hbox{Preprint, arXiv:2405.03626}
\medskip

\title
[On determinants involving Legendre symbols]
{Problems and results on determinants \\ involving Legendre symbols}
\author
[Zhi-Wei Sun] {Zhi-Wei Sun}

\address{School of Mathematics, Nanjing
University, Nanjing 210093, People's Republic of China}
\email{zwsun@nju.edu.cn}

\keywords{Determinants, Legendre symbols, quadratic residues modulo primes.
\newline \indent Supported
by the National Natural Science Foundation of China (grant 12371004).}

\begin{abstract}
In this paper we investigate determinants whose entries are linear combinations of
 Legendre symbols. We deduce some new results in this direction; for example, we prove that
 for any prime $p\equiv3\pmod4$ we have
 $$\det\left[x+\left(\frac{j-k}p\right)+\left(\frac jp\right)-\left(\frac kp\right)\right]_{0\ls j,k\ls(p-1)/2}=4,$$
 where $(\frac{\cdot}p)$ is the Legendre symbol.
 We also pose many conjectures for further research. For example, for any prime $p>3$ we conjecture that
\begin{align*}&\ \det\left[\left(\frac{j+k}p\right)+\left(\frac{j-k}p\right)+\left(\frac{jk}p\right)\right]_{1\ls j,k\ls(p-1)/2}
\\=&\ \begin{cases}(\frac 2p)p^{(p-5)/4}&\text{if}\ p\equiv1\pmod4,
\\(-1)^{(h(-p)-1)/2}(1-(2-(\frac 2p))h(-p))p^{(p-3)/4}&\text{if}\ p\equiv3\pmod4,
\end{cases}\end{align*}
where $h(-p)$ is the class number of the imaginary quadratic field $\mathbb Q(\sqrt{-p})$.
\end{abstract}
\maketitle

\section{Introduction}
\setcounter{lemma}{0}
\setcounter{theorem}{0}
\setcounter{corollary}{0}
\setcounter{remark}{0}
\setcounter{equation}{0}

Let $p$ be an odd prime, and let $(\f{\cdot}p)$ denote the Legendre symbol.
For any integer $a\not\eq0\pmod p$, by the quadratic Gauss sum formula we have
$$\sum_{k=0}^{p-1}e^{2\pi iak^2/p}=\l(\f ap\r)\sqrt{(-1)^{(p-1)/2}p}.$$
Let $\ve_p$ and $h(p)$ be the fundamental unit and the class number of the real quadratic field $\Q(\sqrt p)$.
When $p\eq1\pmod4$, by Dirichlet's class number formula we have
$$\prod_{m=1}^{p-1}(1-e^{2\pi i m/p})^{(\f mp)}=\ve_p^{-2h(p)},$$
which implies that
$$\prod_{k=1}^{(p-1)/2}\l(1-e^{2\pi i ak^2/p}\r)=\sqrt p\,\ve_p^{-(\f ap)h(p)}$$
for each integer $a\not\eq0\pmod p$ (see, e.g., \cite[Theorem 1.3(i)]{S19p}).
For convenience, we write
\begin{equation}
\label{a+b} \ve_p^{h(p)}=a_p+b_p\sqrt p\ \ \t{with}\ 2a_p,2b_p\in\Z.
\end{equation}

For a matrix $A=[a_{jk}]_{1\ls j,k\ls n}$ over a field, let $\det A$ or $|a_{jk}|_{1\ls j,k\ls n}$
denote its determinant. In this paper we focus on determinants involving Legendre symbols.

Let $p=2n+1$ be an odd prime.
In 2004, R. Chapman \cite{Ch} used quadratic
Gauss sums and Dirichlet's class number formula to determine the determinants of the matrices
\begin{equation*}C_p(x)=\l[x+\l(\f{j+k-1}p\r)\r]_{1\ls j,k\ls n}\ \t{and}\
C_p^*(x)=\l[x+\l(\f{j+k-1}p\r)\r]_{1\ls j,k\ls n+1}.
\end{equation*}
By \cite[Corollary 3]{Ch}, provided $p>3$ we have
\begin{equation}\label{C}\det C_p(x)=\begin{cases}(-1)^{n/2}2^n(b_p-a_px)&\t{if}\ p\eq1\pmod4,
\\-2^nx&\t{if}\ p\eq3\pmod4,\end{cases}
\end{equation}
and
\begin{equation}\label{C*}\det C_p^*(x)=\begin{cases}(-1)^{n/2}2^n(pb_px-a_p)&\t{if}\ p\eq1\pmod4,
\\2^n&\t{if}\ p\eq3\pmod4.\end{cases}
\end{equation}
Since $(n+1-j)+(n+1-k)-1\eq-j-k\pmod p$, we also have
\begin{equation}\label{C+}\det C_p(x)=\bg|x+\l(\f{-j-k}p\r)\bg|_{1\ls j,k \ls n}
=(-1)^n\bg|(-1)^nx+\l(\f{j+k}p\r)\bg|_{1\ls j,k\ls n}
\end{equation}
and
\begin{equation}\label{C*+}\det C_p^*(x)=\bg|x+\l(\f{-j-k}p\r)\bg|_{0\ls j,k\ls n}=\bg|(-1)^n x+\l(\f{j+k}p\r)\bg|_{0\ls j,k\ls n}.
\end{equation}

Let $p$ be an odd prime, and write
\begin{equation}\label{ab'}\ve_p^{(2-(\f2p))h(p)}=a_p'+b_p'\sqrt p
\ \ \t{with}\ 2a_p',2b_p'\in\Z.
\end{equation}
In 2003, Chapman conjectured that
$$\bg|\l(\f{j-k}p\r)\bg|_{0\ls j,k\ls(p-1)/2}=\begin{cases}-a_p'&\t{if}\ p\eq1\pmod4,
\\1&\t{if}\ p\eq3\pmod4;\end{cases}$$
this challenging conjecture was finally confirmed by M. Vsemirnov \cite{V12,V13} in 2012-2013 via matrix decomposition
and quadratic Gauss sums.
Recently, L.-Y. Wang, H.-L. Wu and H.-X. Ni \cite{WWN} extended this as follows:
\begin{equation}\label{evil-x}\bg|x+\l(\f{j-k}p\r)\bg|_{0\ls j,k\ls(p-1)/2}=\begin{cases}(\f 2p)pb_p'x-a_p'&\t{if}\ p\eq1\pmod4,
\\1&\t{if}\ p\eq3\pmod4,\end{cases}
\end{equation}
which was ever conjectured by the author.

Let $p=2n+1$ be an odd prime, and let $d\in\Z$. The author \cite{S19} initiated the study of the determinants
$$S(d,p)=\bg|\l(\f{j^2+dk^2}p\r)\bg|_{1\ls j,k\ls n}\ \t{and}\ T(d,p)=\bg|\l(\f{j^2+dk^2}p\r)\bg|_{0\ls j,k\ls n}.$$
He proved that
\begin{equation*}\label{ST} S(d,p)=\begin{cases}\f2{p-1}T(d,p)&\t{if}\ (\f dp)=1,\\0&\t{if}\ (\f dp)=-1.
\end{cases},
\end{equation*}
and
\begin{equation*}\l(\f{T(d,p)}p\r)=\begin{cases}(\f{2}p)&\t{if}\ (\f dp)=1,
\\1&\t{if}\ (\f dp)=-1.\end{cases}
\end{equation*}

We first state a basic result.

\begin{theorem} \label{Th1.1} {\rm (i)} Let $p$ be an odd prime, and let $m,n\in\Z$ with $n\gs m+3$. Then, for any complex numbers $a,b,c,d$, we have
\begin{equation}\label{abcd}\l|a+b\l(\f jp\r)+c\l(\f kp\r)+d\l(\f{jk}p\r)\r|_{m\ls j,k\ls n}=0.
\end{equation}

{\rm (ii)} Let $p>5$ be a prime with $p\eq1\pmod4$. For any $\da\in\{\pm1\}$ and $m\in\{0,1\}$, we have
\begin{equation}\label{j^2+-k^2}\bg|x+\l(\f{j^2+k^2}p\r)+\da\l(\f{j^2-k^2}p\r)\bg|_{m\ls j,k\ls (p-1)/2}=0.
\end{equation}
\end{theorem}
\begin{remark} In 1956, D. H. Lehmer \cite{L} found all the eigenvalues of the determinant
$$\l|a+b\l(\f jp\r)+c\l(\f kp\r)+d\l(\f{jk}p\r)\r|_{1\ls j,k\ls p-1},$$
where $p$ is an odd prime and $a,b,c,d$ are complex numbers.
As a supplement to Theorem \ref{Th1.1}(ii), we conjecture that
\begin{equation}\l|x+\l(\f{j^2+k^2}p\r)+\l(\f{j^2-k^2}p\r)\r|_{1\ls j,k\ls (p-1)/2}
=\l(\f{p-1}2x-1\r)p^{(p-3)/4}
\end{equation}
for any prime $p\eq3\pmod4$.
\end{remark}

Now we state our central result.

\begin{theorem} \label{Th1.2} Let $p$ be an odd prime, and let $a_i,b_i\in\Z$ for all $i=1,\ldots,m$.
Let $c_1,\ldots,c_m$ be complex numbers, and set
$$c=\sum_{s=1}^m c_s\l(\f{a_s}p\r)\sum_{t=1}^m c_t\l(\f{b_t}p\r).$$

{\rm (i)} For each $n\in\{1,\ldots,p-1\}$, we have
\begin{equation} \label{cx} \l|\sum_{i=1}^mc_i\l(\f{a_ij+b_ik}p\r)+\l(\f{jk}p\r)x\r|_{0\ls j,k\ls n}
=\l|\sum_{i=1}^mc_i\l(\f{a_ij+b_ik}p\r)\r|_{0\ls j,k\ls n}
\end{equation}
and
\begin{equation}\label{01}\begin{aligned}&\ c\l|\sum_{i=1}^mc_i\l(\f{a_ij+b_ik}p\r)+\l(\f{jk}p\r)x\r|_{1\ls j,k\ls n}
\\=&\ c\l|\sum_{i=1}^mc_i\l(\f{a_ij+b_ik}p\r)\r|_{1\ls j,k\ls n}
- x\l|\sum_{i=1}^mc_i\l(\f{a_ij+b_ik}p\r)\r|_{0\ls j,k\ls n}.
\end{aligned}
\end{equation}

{\rm (ii)} For any positive integer $n$, we have
\begin{equation}\label{yz}\begin{aligned}&\ c\l|x+\sum_{i=1}^m c_i\l(\f{a_ij+b_ik}p\r)+\l(\f jp\r)y+\l(\f kp\r)z\r|_{0\ls j,k\ls n}
\\=&\ \l( y+\sum_{i=1}^m c_i\l(\f{a_i}p\r)\r)\l(z+\sum_{i=1}^m c_i\l(\f{b_i}p\r)\r)\l|\sum_{i=1}^m c_i\l(\f{a_ij+b_ik}p\r)\r|_{0\ls j,k\ls n}
\\&\ +cx\l|\sum_{i=1}^m c_i\l(\f{a_ij+b_ik}p\r)-\sum_{i=1}^m c_i\l(\f{a_i}p\r)\l(\f jp\r)-\sum_{i=1}^m c_i\l(\f{b_i}p\r)\l(\f kp\r)\r|_{1\ls j,k\ls n}.
\end{aligned}
\end{equation}
\end{theorem}

Applying Theorem \ref{Th1.2} and using the known values of
$$\l|x+\l(\f{j+k}p\r)\r|_{0\ls j,k\ls(p-1)/2}\ \t{and}\ \l|x+\l(\f{j-k}p\r)\r|_{0\ls j,k\ls(p-1)/2}$$
(where $p$ is an odd prime), we can deduce the following general result.

\begin{theorem}\label{Th1.3} Let $p$ be an odd prime.

{\rm (i)} If $p>3$, then
\begin{equation}\label{+}\begin{aligned}&\ \l|x+\l(\f{j+k}p\r)+\l(\f jp\r)y+\l(\f kp\r)z\r|_{0\ls j,k\ls(p-1)/2}
\\=&\ \begin{cases}(\f 2p)2^{(p-1)/2}(pb_px-(y+1)(z+1)a_p)&\t{if}\ p\eq1\pmod4,
\\(y+1)(z+1)2^{(p-1)/2}&\t{if}\ p\eq3\pmod4.
\end{cases}\end{aligned}
\end{equation}

{\rm (ii)} We have
\begin{equation}\label{-}\begin{aligned}&\ \l|x+\l(\f{j-k}p\r)+\l(\f jp\r)y+\l(\f kp\r)z\r|_{0\ls j,k\ls(p-1)/2}
\\=&\ \begin{cases}(\f 2p)pb_p'x-(1+y)(1+z)a_p'&\t{if}\ p\eq1\pmod 4,
\\(1+y)(1-z)&\t{if}\ p\eq3\pmod4.\end{cases}
\end{aligned}
\end{equation}
\end{theorem}

Clearly, this theorem has the following consequence.

\begin{corollary} Let $p$ be a prime with $p\eq3\pmod4$. Then
\begin{equation}\l|x+\l(\f{j-k}p\r)+\l(\f jp\r)-\l(\f kp\r)\r|_{0\ls j,k\ls(p-1)/2}=4.
\end{equation}
When $p>3$, we have
\begin{equation}\l|x+\l(\f{j+k}p\r)+\l(\f jp\r)+\l(\f kp\r)\r|_{0\ls j,k\ls(p-1)/2}=2^{(p+3)/2}.
\end{equation}
\end{corollary}

We also have the following general result.

\begin{theorem}\label{Th1.4} Let $p>3$ be a prime. Then
\begin{equation}\label{1.4-1}\begin{aligned}&\ \l|x+\l(\f{j+k-1}p\r)+\l(\f jp\r)y+\l(\f kp\r)z\r|_{1\ls j,k\ls(p-1)/2}
\\=&\ \begin{cases}(\f 2p)2^{(p-1)/2}((yz-x)a_p+(y+1)(z+1)b_p)&\t{if}\ p\eq1\pmod4,
\\2^{(p-1)/2}(yz-x)&\t{if}\ p\eq3\pmod4.
\end{cases}
\end{aligned}\end{equation}
Also,
\begin{equation}\label{1.4-2}\begin{aligned}&\ \l|x+\l(\f{j+k-1}p\r)+\l(\f jp\r)y+\l(\f kp\r)z\r|_{1\ls j,k\ls(p+1)/2}
\\=&\ \begin{cases}(\f 2p)2^{(p-1)/2}(pb_p(x-yz)-a_p(y+1)(z+1))&\t{if}\ p\eq1\pmod4,
\\2^{(p-1)/2}(y+1)(z+1)&\t{if}\ p\eq3\pmod4.\end{cases}
\end{aligned}\end{equation}
\end{theorem}

We are going to prove Theorems 1.1-1.4 in the next section, and pose in Sections 3-5 many conjectures on determinants involving linear combinations of Legendre symbols.

\section{Proofs of Theorems \ref{Th1.1}-\ref{Th1.3}}
\setcounter{lemma}{0}
\setcounter{theorem}{0}
\setcounter{corollary}{0}
\setcounter{remark}{0}
\setcounter{equation}{0}

\medskip
\noindent{\tt Proof of Theorem \ref{Th1.1}}. (i) We now prove part (i) of Theorem \ref{Th1.1}.
As the four Legendre symbols
$$\l(\f mp\r),\ \l(\f{m+1}p\r),\ \l(\f{m+2}p\r),\ \l(\f{m+3}p\r)$$
cannot be pairwise distinct, there are $j,j'\in\{m,m+1,m+2,m+3\}$ with $j\not=j'$
such that $(\f jp)=(\f{j'}p)$. Thus
$$a+b\l(\f jp\r)+c\l(\f kp\r)+d\l(\f{jk}p\r)=a+b\l(\f {j'}p\r)+c\l(\f kp\r)+d\l(\f{{j'}k}p\r)$$
for all $k=m,\ldots,n$, and hence \eqref{abcd} holds.

(ii) We now turn to prove part (ii) of Theorem \ref{Th1.1}.
Set $n=(p-1)/2$ and $q=n!$. By Wilson's theorem,
$$-1\eq (p-1)!=\prod_{k=1}^n k(p-k)\eq(-1)^n (n!)^2=q^2\pmod p.$$
For each $k=1,\ldots,n$, there is a unique $r_k\in\{1,\ldots,n\}$ such that $qk$ is congruent to $r_k$ or $-r_k$ modulo $p$. Note that $r_k^2\eq-k^2\pmod p$ and $r_k\not=k$ since $q^2\eq-1\not\eq1\pmod p$.
As $qr_k\eq \pm q^2k\eq\mp k\pmod p$, we also have $r_{r_k}=k$. For any $k\in\{1,\ldots,n\}$
and $j\in\{m,\ldots,n\}$, clearly
\begin{align*} &\ x+\l(\f{j^2+k^2}p\r)+\da\l(\f{j^2-k^2}p\r)-\da\l(x+\l(\f{j^2+r_k^2}p\r)+\da\l(\f{j^2-r_k^2}p\r)\r)
\\=&\ x+\l(\f{j^2+k^2}p\r)+\da\l(\f{j^2+r_k^2}p\r)
-\da\l(x+\l(\f{j^2+r_k^2}p\r)+\da\l(\f{j^2+k^2}p\r)\r)
\\=&\ (1-\da)x.
\end{align*}
When $\da=1$, this clearly implies the equality \eqref{j^2+-k^2}.

Now we consider the case $\da=-1$. As $n=(p-1)/2\gs 4$, we may choose $k\in\{1,\ldots,n\}\sm\{1,r_1\}$.
Note that
$1,r_1,k,r_k$ are distinct elements of $\{1,\ldots,n\}$ with
\begin{align*}&\ x+\l(\f{j^2+k^2}p\r)+\da\l(\f{j^2-k^2}p\r)-\da\l(x+\l(\f{j^2+r_k^2}p\r)+\da\l(\f{j^2-r_k^2}p\r)\r)
\\=&\ 2x=x+\l(\f{j^2+1^2}p\r)+\da\l(\f{j^2-1^2}p\r)-\da\l(x+\l(\f{j^2+r_1^2}p\r)+\da\l(\f{j^2-r_1^2}p\r)\r)
\end{align*}
for all $j=m,\ldots,n$. So \eqref{j^2+-k^2} holds.

In view of the above, we have completed the proof of Theorem \ref{Th1.1}. \qed

To prove Theorem \ref{Th1.2}, we need the following basic lemma which can be found in \cite[Lemma 2.1]{RJ}.

\begin{lemma} \label{Lem2.1} Let $A=[a_{jk}]_{0\ls j,k\ls m}$ be a matrix over a field. Then
\begin{equation}\label{ab}\det[x+a_{jk}]_{0\ls j,k\ls m}-\det[a_{jk}]_{0\ls j,k\ls m}
=x\det[b_{jk}]_{1\ls j,k\ls m},
\end{equation}
where $b_{jk}=a_{jk}-a_{j0}-a_{0k}+a_{00}$.
\end{lemma}

\medskip
\noindent{\tt Proof of Theorem \ref{Th1.2}}. For convenience, we set
$$f(j,k)=\sum_{i=1}^m c_i\l(\f{a_ij+b_ik}p\r)$$
for any $j,k=0,1,2,\ldots$.

(i) We first prove part (i) of Theorem \ref{Th1.2}.
For $i=1,\ldots,m$ and $j,k=1,\ldots,n$, clearly
$$f(j,k)+\l(\f{jk}p\r)x=\l(\f jp\r)\l(\f kp\r)\l(\l(\f{jk}p\r)f(j,k)+x\r).$$
It follows that
\begin{equation}\label{x^2} x^2\l|f(j,k)+\l(\f{jk}p\r)x\r|_{0\ls j,k\ls n}
 =c\prod_{j=1}^n\l(\f jp\r)
\times\prod_{k=1}^n\l(\f kp\r)\times\det A_0=c\det A_0,
\end{equation}
where $A_0$ is obtained from the matrix $A=[(\f{jk}p)f(j,k)+x]_{0\ls j,k\ls n}$ via replacing the first entry $x$ in the first row by $0$. If we expand $\det A_0$ and $\det A$ along their first rows, we immediately see that
\begin{equation}\label{A-A0}\begin{aligned}\det A-\det A_0&=x\l|\l(\f{jk}p\r)f(j,k)+x\r|_{1\ls j,k\ls n}.
\end{aligned}\end{equation}
By Lemma \ref{Lem2.1},
$$\det A=\l|\l(\f{jk}p\r)f(j,k)\r|_{0\ls j,k\ls n}+x\l|\l(\f{jk}p\r)f(j,k)\r|_{1\ls j,k\ls n}=x\l|f(j,k)\r|_{1\ls j,k\ls n}.$$
Combining this with \eqref{x^2} and \eqref{A-A0}, we obtain
\begin{align*}&\ x^2\l|f(j,k)+\l(\f{jk}p\r)x\r|_{0\ls j,k\ls n}
\\=&\ c\l(x\l|f(j,k)\r|_{1\ls j,k\ls n}-x\l|\l(\f{jk}p\r)f(j,k)+x\r|_{1\ls j,k\ls n}\r)
\end{align*}
and hence
\begin{equation}\label{f(j,k)}\begin{aligned}&\ x\l|f(j,k)+\l(\f{jk}p\r)x\r|_{0\ls j,k\ls n}
\\=&\ c\l(\l|f(j,k)\r|_{1\ls j,k\ls n}-\l|\l(\f{jk}p\r)f(j,k)+x\r|_{1\ls j,k\ls n}\r).
\end{aligned}
\end{equation}

Applying Lemma \ref{Lem2.1}, we find that
\begin{align*}
\l|\l(\f{jk}p\r)f(j,k)+x\r|_{1\ls j,k\ls n}= \l|\l(\f{jk}p\r)f(j,k)\r|_{1\ls j,k\ls n}+x\det D,
\end{align*}
where $D=[d_{jk}]_{2\ls j,k\ls n}$ with
\begin{align*}d_{jk}=&\ \l(\f{jk}p\r)f(j,k)-\l(\f jp\r)f(j,1)-\l(\f kp\r)f(1,k)+f(1,1).
\end{align*}
Therefore
\begin{align*} \l|\l(\f{jk}p\r)f(j,k)+x\r|_{1\ls j,k\ls n}
=\l|f(j,k)\r|_{1\ls j,k\ls n}+x\det D.
\end{align*}
Combining this with \eqref{f(j,k)}, we immediately get
$$\l|f(j,k)+\l(\f{jk}p\r)x\r|_{0\ls j,k\ls n}=-c\det D$$
and hence \eqref{cx} follows. In light of \eqref{cx} and \eqref{f(j,k)},
$$c|f(j,k)|_{1\ls j,k\ls n}-x|f(j,k)|_{0\ls j,k\ls n}=c\l|f(j,k)+\l(\f{jk}p\r)x\r|_{1\ls j,k\ls n}
,$$ which gives \eqref{01}.

(ii) Now we turn to prove part (ii) of Theorem \ref{Th1.2}. Let
$$a_{jk}=f(j,k)+\l(\f jp\r)y+\l(\f kp\r)z$$
for $j,k=0,\ldots,n$. It is easy to see that
$$a_{jk}-a_{j0}-a_{0k}+a_{00}=f(j,k)-\sum_{i=1}^m c_i\l(\f{a_i}p\r)\l(\f jp\r)
-\sum_{i=1}^m c_i\l(\f{b_i}p\r)\l(\f kp\r).$$
Thus, in view of Lemma \ref{Lem2.1},
\begin{align*}&\ |x+a_{jk}|_{0\ls j,k\ls n}-|a_{jk}|_{0\ls j,k\ls n}
\\=&\ x\l|f(j,k)-\sum_{i=1}^m c_i\l(\f{a_i}p\r)\l(\f jp\r)-\sum_{i=1}^m c_i\l(\f{b_i}p\r)\l(\f kp\r)\r|_{1\ls j,k\ls n}.
\end{align*}
So we have reduced \eqref{yz} to the equality
\begin{equation}\label{no-x}c\l|a_{jk}\r|_{0\ls j,k\ls n}=\l(y+\sum_{i=1}^m c_i\l(\f{a_i}p\r)\r)
\l(z+\sum_{i=1}^m c_i\l(\f{b_i}p\r)\r)|f(j,k)|_{0\ls j,k\ls n}.
\end{equation}

For $k=0,\ldots,n$, clearly
$$a_{0k}=f(0,k)+\l(\f kp\r)z=\l(z+\sum_{i=1}^m c_i\l(\f{b_i}p\r)\r)\l(\f kp\r)$$
and
$$a_{jk}-a_{0k}=f(j,k)+\l(\f jp\r)y-\sum_{i=1}^m c_i\l(\f{b_i}p\r)\l(\f kp\r)$$
for all $j=1,\ldots,n$. Thus
\begin{equation}\label{bi}\sum_{i=1}^m c_i\l(\f{b_i}p\r)\times|a_{jk}|_{0\ls j,k\ls n}
=\l(z+\sum_{i=1}^m c_i\l(\f{b_i}p\r)\r)\l|f(j,k)+\l(\f jp\r)y\r|_{0\ls j,k\ls n}.
\end{equation}
For $j=0,\ldots,n$, apparently
$$f(j,0)+\l(\f jp\r)y=\l(y+\sum_{i=1}^m c_i\l(\f{a_i}p\r)\r)\l(\f jp\r)$$
and
$$f(j,k)+\l(\f jp\r)y-\l(f(j,0)+\l(\f jp\r)y\r)=f(j,k)-\sum_{i=1}^m c_i\l(\f{a_i}p\r)\l(\f jp\r)$$
for all $k=1,\ldots,n$. Therefore
$$\sum_{i=1}^m c_i\l(\f{a_i}p\r)\times\l|f(j,k)+\l(\f jp\r)y\r|_{0\ls j,k\ls n}
=\l(y+\sum_{i=1}^m c_i\l(\f{a_i}p\r)\r)|f(j,k)|_{0\ls j,k\ls n}.$$
Combining this with \eqref{bi}, we immediately obtain the desired \eqref{no-x}.

In view of the above, we have completed the proof of Theorem \ref{Th1.2}.
\qed

Recall that an $n\times n$ matrix $A=[a_{jk}]_{1\ls j,k\ls n}$ over a field
is called {\it skew-symmetric} if $a_{jk}+a_{kj}=0$ for all $j,k=1,\ldots,n$.

Suppose that $A=[a_{jk}]_{1\ls j,k\ls n}$ is a skew-symmetric matrix over $\Z$. Note that
$$\det A=|a_{kj}|_{1\ls j,k\ls n}=|-a_{jk}|_{1\ls j,k\ls n}=(-1)^n\det A.$$
Thus $\det A=0$ if $n$ is odd. By a theorem of Cayley, $\det A$ is a square if $n$ is even
(cf. \cite{K05}).

\begin{lemma}\label{Lem2.2} Let $p$ be an odd prime. Then
\begin{equation}\label{2.2-1}\begin{aligned}&\ \l|\l(\f{j+k}p\r)-\l(\f jp\r)-\l(\f kp\r)\r|_{1\ls j,k\ls(p-1)/2}
\\=&\ \begin{cases}(\f 2p)2^{(p-1)/2}pb_p&\t{if}\ p\eq1\pmod4,
\\0&\t{if}\ p>3\ \t{and}\ p\eq3\pmod4.\end{cases}
\end{aligned}
\end{equation}
We also have
\begin{equation}\label{2.2-2}\l|\l(\f{j-k}p\r)-\l(\f jp\r)-\l(\f {-k}p\r)\r|_{1\ls j,k\ls(p-1)/2}
=\begin{cases}(\f 2p)pb_p'&\t{if}\ p\eq1\pmod4,
\\0&\t{if}\ p\eq3\pmod4.\end{cases}
\end{equation}
\end{lemma}
\Proof. Let $n=(p-1)/2$ and $\da\in\{\pm1\}$.
Define $a_{jk}=(\f{j+\da k}p)$ for $j,k=0,\ldots,n$. Then
$$a_{jk}-a_{j0}-a_{0k}+a_{00}=\l(\f{j+\da k}p\r)-\l(\f jp\r)-\l(\f{\da k}p\r).$$
Thus, by Lemma \ref{Lem2.1} we have
\begin{equation}\label{a-delta}\begin{aligned}&\ \det[1+a_{jk}]_{0\ls j,k\ls n}-\det[a_{jk}]_{1\ls j,k\ls n}
\\=&\ \l|\l(\f{j+\da k}p\r)-\l(\f jp\r)-\l(\f{\da k}p\r)\r|_{1\ls j,k\ls n}.
\end{aligned}\end{equation}

Combining \eqref{C*} and \eqref{C*+}, we obtain
\begin{equation}\label{2.10}\l|x+\l(\f{j+k}p\r)\r|_{0\ls j,k\ls n}
=\begin{cases}(\f 2p)2^n(pb_px-a_p)&\t{if}\ p\eq1\pmod4,
\\2^n&\t{if}\ p>3\ \&\ p\eq3\pmod4.\end{cases}
\end{equation}
So we know the exact value of $|x+a_{jk}|_{0\ls j,k\ls n}$ in the case $\da=1$.
When $\da=-1$, the equality \eqref{evil-x} gives the exact value of $|x+a_{jk}|_{0\ls j,k\ls n}$.
Since $|x+a_{jk}|_{0\ls j,k\ls n}$ is evaluated, we immediately obtain the exact value of
$$\l|\l(\f{j+\da k}p\r)-\l(\f jp\r)-\l(\f{\da k}p\r)\r|_{1\ls j,k\ls n}$$
by using \eqref{a-delta}. Therefore \eqref{2.2-1} and \eqref{2.2-2} hold.
In the case $p\eq3\pmod4$, we may prove \eqref{2.2-2} without using \eqref{evil-x} since
the matrix in \eqref{2.2-2} is skew-symmetric and of odd order. This ends our proof. \qed

\medskip
\noindent{\tt Proof of Theorem \ref{Th1.3}}. Set $n=(p-1)/2$. Combining Theorem \ref{Th1.2}(ii),
\eqref{2.10}, \eqref{evil-x} and Lemma \ref{Lem2.2}, we immediately obtain the desired results.
\qed

\medskip
\noindent{\tt Proof of Theorem \ref{Th1.4}}. Let $n\in\{(p-1)/2,(p+1)/2\}$, and set
$$a_{jk}=x+\l(\f{j+k-1}p\r)+\l(\f jp\r)y$$
for $j,k=1,\ldots,n$. Observe that
$$a_{jk}-a_{j1}+\f{a_{j1}}{y+1}=\l(\f {j+k-1}p\r)+\f x{y+1}$$
for all $1<j\ls n$ and $1<k\ls n$. Thus
\begin{align*}&\ \l|x+\l(\f{j+k-1}p\r)+\l(\f jp\r)y\r|_{1\ls j,k\ls n}
\\=&\ (y+1)\l|\l(\f{j+k-1}p\r)+\f x{y+1}\r|_{1\ls j,k\ls n}.
\end{align*}
Combining this with \eqref{C} and \eqref{C*}, we have
\begin{equation}\label{1y}\begin{aligned}&\ \l|x+\l(\f{j+k-1}p\r)+\l(\f jp\r)y\r|_{1\ls j,k\ls (p-1)/2}
\\=&\ \begin{cases}(\f 2p)2^{(p-1)/2}(b_p(y+1)-a_px)&\t{if}\ p\eq1\pmod4,
\\-2^{(p-1)/2}x&\t{if}\ p\eq3\pmod4,\end{cases}
\end{aligned}
\end{equation}
and
\begin{equation}\label{2y}\begin{aligned}&\ \l|x+\l(\f{j+k-1}p\r)+\l(\f jp\r)y\r|_{1\ls j,k\ls (p+1)/2}
\\=&\ \begin{cases}(\f 2p)2^{(p-1)/2}(pb_px-a_p(y+1))&\t{if}\ p\eq1\pmod4,
\\2^{(p-1)/2}(y+1)&\t{if}\ p\eq3\pmod4.\end{cases}
\end{aligned}
\end{equation}

Let
$$b_{jk}=x+\l(\f{j+k-1}p\r)+\l(\f jp\r)y+\l(\f kp\r)z$$
for $j,k=1,\ldots,n$. Note that
$$b_{jk}-b_{1k}+\f{b_{1k}}{z+1}=\l(\f{j+k-1}p\r)+\l(\f jp\r)y+\f{x-yz}{z+1}$$
for all $j,k=1,\ldots,n$.
Thus
$$|b_{jk}|_{1\ls j,k\ls n}=(z+1)\l|\l(\f{j+k-1}p\r)+\l(\f jp\r)y+\f{x-yz}{z+1}\r|_{1\ls j,k\ls n}.$$
Combining this with \eqref{1y} and \eqref{2y}, we immediately obtain the desired identities
\eqref{1.4-1} and \eqref{1.4-2}. Note that both sides of the equalities \eqref{1.4-1} and \eqref{1.4-2}
 are polynomials in $x,y,z$. If we view $y$ and $z$ as complex numbers, to handle the case $y=-1$ or $z=-1$ we may take limits. This concludes our proof of Theorem \ref{Th1.4}. \qed

\newpage
\section{Conjectures on determinants involving
 $(\f{j\pm k}p)$, $(\f jp)$, $(\f kp)$ and $(\f{jk}p)$}
\setcounter{lemma}{0}
\setcounter{theorem}{0}
\setcounter{corollary}{0}
\setcounter{remark}{0}
\setcounter{conjecture}{0}
\setcounter{equation}{0}

\begin{conjecture}\label{Conj0(p-1)}
Let $p$ be an odd prime.

{\rm (i)} If $p>3$, then
\begin{equation}\begin{aligned}&\ \l|x+\l(\f{j+k}p\r)+\l(\f jp\r)y+\l(\f kp\r)z+\l(\f{jk}p\r)w\r|_{0\ls j,k\ls(p-1)/2}
\\=&\ \begin{cases}(\f 2p)2^{(p-1)/2}(pb_px+a_p(wx-(y+1)(z+1)))&\t{if}\ p\eq1\pmod4,
\\2^{(p-1)/2}((y+1)(z+1)-wx)&\t{if}\ p\eq3\pmod4.\end{cases}
\end{aligned}
\end{equation}

{\rm (ii)} We have
\begin{equation}\begin{aligned}&\ \l|x+\l(\f{j-k}p\r)+\l(\f jp\r)y+\l(\f kp\r)z+\l(\f{jk}p\r)w\r|_{0\ls j,k\ls(p-1)/2}
\\=&\ \begin{cases}a_p'(wx-(y+1)(z+1))+(\f 2p)pb_p'x&\t{if}\ p\eq1\pmod4,
\\wx+(1+y)(1-z)&\t{if}\ p\eq3\pmod4.\end{cases}
\end{aligned}
\end{equation}
\end{conjecture}
\begin{remark} Conjecture \ref{Conj0(p-1)} in the case $wx=0$ follows from Theorems \ref{Th1.2} and \ref{Th1.3}.
\end{remark}

\begin{conjecture}\label{Conj0(p-3)}
Let $p$ be an odd prime, and set $v=wx-(y+1)(z+1)$.

{\rm (i)} If $p>3$, then
\begin{equation}\begin{aligned}&\ \l|x+\l(\f{j+k}p\r)+\l(\f jp\r)y+\l(\f kp\r)z+\l(\f{jk}p\r)w\r|_{0\ls j,k\ls(p-3)/2}
\\=&\ \begin{cases}(\f 2p)2^{(p-3)/2}((pb_p-2a_p)x+(a_p-2b_p)v)&\t{if}\ p\eq1\pmod4,
\\2^{(p-3)/2}(v-2x)&\t{if}\ p\eq3\pmod4.\end{cases}
\end{aligned}
\end{equation}

{\rm (ii)} We have
\begin{equation}\begin{aligned}&\ \l|x+\l(\f{j-k}p\r)+\l(\f jp\r)y+\l(\f kp\r)z+\l(\f{jk}p\r)w\r|_{0\ls j,k\ls(p-3)/2}
\\=&\ \begin{cases}-a_p'x-(\f 2p)b_p'v&\t{if}\ p\eq1\pmod4,
\\x&\t{if}\ p\eq3\pmod4.\end{cases}
\end{aligned}
\end{equation}
\end{conjecture}
\begin{remark} In light of Theorem 1.2, in the case $wx=0$ we can reduce Conjecture \ref{Conj0(p-3)} to the case $y=z=0$.
For any prime $p\eq3\pmod4$, clearly $|(\f{j-k}p)|_{1\ls j,k\ls(p-1)/2}=0$ since the matrix is skew-symmetric and of odd order; the author \cite{S19} conjectured that
$$\l|x+\l(\f{j-k}p\r)\r|_{1\ls j,k\ls(p-1)/2}=x,
\ \ \t{i.e.,}\ \l|x+\l(\f{j-k}p\r)\r|_{0\ls j,k\ls(p-3)/2}=x.$$
In view of Lemma \ref{Lem2.1} or Theorem 1.2(ii), for any prime $p>3$, part (i) of Conjecture \ref{Conj0(p-3)} implies that
\begin{equation}\label{+--}\begin{aligned}&\ \l|\l(\f{j+k}p\r)-\l(\f jp\r)-\l(\f kp\r)\r|_{1\ls j,k\ls(p-3)/2}
\\=&\ \begin{cases}(\f 2p)2^{(p-3)/2}(pb_p-2a_p)&\t{if}\ p\eq1\pmod4,
\\-2^{(p-1)/2}&\t{if}\ p\eq3\pmod4,\end{cases}
\end{aligned}\end{equation}
while
part (ii) of Conjecture \ref{Conj0(p-3)} implies that
\begin{equation}\label{S0}\l|\l(\f{j-k}p\r)-\l(\f jp\r)-\l(\f {-k}p\r)\r|_{1\ls j,k\ls(p-3)/2}=
\begin{cases}-a_p'&\t{if}\ p\eq1\pmod4,\\1&\t{if}\ p\eq3\pmod4.\end{cases}
\end{equation}
\end{remark}

\begin{conjecture}\label{Conj-(p-5)} Let $p\gs 5$ be a prime. Then
\begin{equation}\begin{aligned}&\l|x+\l(\f{j-k}p\r)+\l(\f jp\r)y+\l(\f kp\r)z+\l(\f{jk}p\r)w\r|_{0\ls j,k\ls(p-5)/2}
\\=&\ \begin{cases}(\f 2p)(2a_p'-pb_p')x+(a_p'-2b_p')((1+y)(1+z)-wx)&\t{if}\ p\eq1\pmod4,
\\wx+(1+y)(1-z)&\t{if}\ p\eq3\pmod4.\end{cases}
\end{aligned}
\end{equation}
\end{conjecture}
\begin{remark} In light of Theorem 1.2, in the case $wx=0$ we can reduce this conjecture to the case $y=z=0$.
In a unpublished preprint written in 2003, for each prime $p\gs5$ with $p\eq3\pmod4$, R. Chapman conjectured that
$$\l|\l(\f{j-k}p\r)\r|_{1\ls j,k\ls(p-3)/2}=1,
 \ \t{i.e.},\ \l|\l(\f{j-k}p\r)\r|_{0\ls j,k\ls(p-5)/2}=1.$$
In view of Lemma \ref{Lem2.1} of Theorem \ref{Th1.2}(ii), Conjecture \ref{Conj-(p-5)} implies that
\begin{equation}\l|\l(\f{j-k}p\r)-\l(\f jp\r)-\l(\f {k}p\r)\r|_{1\ls j,k\ls(p-5)/2}=\l(\f 2p\r)(2a_p'-pb_p')
\end{equation}
for any prime $p>5$ with $p\eq1\pmod4$, and that
\begin{equation}\label{0(p-5)}\l|\l(\f{j-k}p\r)-\l(\f jp\r)+\l(\f {k}p\r)\r|_{1\ls j,k\ls(p-5)/2}=0
\end{equation}
for any prime $p>5$ with $p\eq3\pmod4$. The equality \eqref{0(p-5)} is easy since
 the matrix is skew-symmetric and of odd order.
\end{remark}

\begin{conjecture} For any prime $p\gs7$ with $p\eq3\pmod4$, we have
\begin{equation}\l|x+\l(\f{j-k}p\r)\r|_{0\ls j,k\ls(p-7)/2}=\l\lfloor\f{p-2}3\r\rfloor^2 x.
\end{equation}
\end{conjecture}
\begin{remark} Surprisingly, this concise conjecture has not been found before.
\end{remark}

\begin{conjecture} Let $p>3$ be a prime.

{\rm (i)} If $p\eq1\pmod4$, then
\begin{equation}\begin{aligned}&\ \l|x+\l(\f{j+k}p\r)+\l(\f jp\r)y +\l(\f kp\r)z+\l(\f{jk}p\r)w\r|_{1\ls j,k\ls(p-1)/2}
\\=&\ \l(\f 2p\r)2^{(p-1)/2}(a_p(w-x)+b_p+(b_p-1)(y+z)-c_p(wx-yz)),
\end{aligned}
\end{equation}
where $c_p=(p+1)b_p-2$.
When $p\eq3\pmod4$, we have
\begin{equation}\begin{aligned}&\ \l|x+\l(\f{j+k}p\r)+\l(\f jp\r)y +\l(\f kp\r)z+\l(\f{jk}p\r)w\r|_{1\ls j,k\ls(p-1)/2}
\\=&\ -2^{(p-1)/2}(w+x+(-1)^{(h(-p)-1)/2}(y+z+2yz-2wx)),
\end{aligned}
\end{equation}
where $h(-p)$ denotes the class number of the imaginary quadratic field $\Q(\sqrt{-p})$.

{\rm (ii)} If $p\eq1\pmod4$, then
\begin{equation}\begin{aligned}&\ \l|x+\l(\f{j-k}p\r)+\l(\f jp\r)y +\l(\f kp\r)z+\l(\f{jk}p\r)w\r|_{1\ls j,k\ls(p-1)/2}
\\=&\ a_p'(w-x)+\l(\f 2p\r)(b_p'+(b_p'-1)(y+z)+c_p'(yz-wx)),
\end{aligned}
\end{equation}
where $c_p'=(p+1)b_p'-2$. When $p\eq3\pmod4$, we have
\begin{equation}\begin{aligned}&\ \l|x+\l(\f{j-k}p\r)+\l(\f jp\r)y +\l(\f kp\r)z+\l(\f{jk}p\r)w\r|_{1\ls j,k\ls(p-1)/2}
\\=&\ w+x-(-1)^{(h(-p)-1)/2}(y+z).
\end{aligned}
\end{equation}
\end{conjecture}

\begin{conjecture} Let $p\gs5$ be a prime.

{\rm (i)} If $p\eq1\pmod4$, then
\begin{equation}\begin{aligned}&\ \l|x+\l(\f{j+k}p\r)+\l(\f jp\r)y+\l(\f kp\r)z+\l(\f{jk}p\r)w\r|_{1\ls j,k\ls(p-3)/2}
\\=&\ \l(\f 2p\r)2^{(p-3)/2}(b_p-a_px+(a_p-2b_p)w+(b_p-1)(y+z)+d_p(yz-wx)),
\end{aligned}
\end{equation}
where $d_p=(p+1)b_p-2(a_p+1)$. When $p\eq3\pmod4$, we have
\begin{equation}\begin{aligned}&\ \l|x+\l(\f{j+k}p\r)+\l(\f jp\r)y+\l(\f kp\r)z+\l(\f{jk}p\r)w\r|_{1\ls j,k\ls(p-3)/2}
\\=&\ 2^{(p-3)/2}(w+x+2(wx-yz)+(-1)^{(h(-p)-1)/2}(y+z+2yz-2wx)).
\end{aligned}
\end{equation}

{\rm (ii)} If $p\eq1\pmod4$, then
\begin{equation}\begin{aligned}&\ \l|x+\l(\f{j-k}p\r)+\l(\f jp\r)y+\l(\f kp\r)z+\l(\f{jk}p\r)w\r|_{1\ls j,k\ls(p-3)/2}
\\=&\ e_p'+\l(\f 2p\r)((2a_p'-pb_p')x-b_p'w)+(e_p'+1)(y+z)+2(b_p'-1)(wx-yz),
\end{aligned}
\end{equation}
where $e_p'=a_p'-2b_p'$. When $p\eq3\pmod4$, we have
\begin{equation}\begin{aligned}&\ \l|x+\l(\f{j-k}p\r)+\l(\f jp\r)y+\l(\f kp\r)z+\l(\f{jk}p\r)w\r|_{1\ls j,k\ls(p-3)/2}
\\=&\ 1+\l(1-(-1)^{(h(-p)-1)/2}\l(\f 2p\r)\r)(2(wx-yz)+y-z).
\end{aligned}
\end{equation}
\end{conjecture}

\begin{conjecture} Let $p>3$ be a prime.

{\rm (i)} We have
\begin{equation}\begin{aligned}&\ \l|x+\l(\f{j+k}p\r)+\l(\f{j-k}p\r)+\l(\f jp\r)y+\l(\f kp\r)z+\l(\f{jk}p\r)w\r|_{0\ls j,k\ls (p-1)/2}
\\=&\ \begin{cases}(\f 2p)p^{(p+3)/4}x&\t{if}\ p\eq1\pmod4,
\\(-1)^{(h(-p)-1)/2}p^{(p-3)/4}(px+(2-(\f 2p))h(-p)v)&\t{if}\ p\eq3\pmod4,
\end{cases}
\end{aligned}
\end{equation}
where $v=(y+2)z-wx$.

{\rm (ii)} If $p\eq1\pmod4$, then
\begin{equation} \l|\l(\f{j+k}p\r)-\l(\f{j-k}p\r)+\l(\f jp\r)y+\l(\f kp\r)z\r|_{0\ls j,k\ls(p-1)/2}=4p^{(p-5)/4}x_pyz
\end{equation}
for some $x_p\in\Z$ only depending on $p$.
\end{conjecture}
\begin{remark} Our computation indicates that
$$x_5=1,\, x_{13}=-3,\, x_{17}=2,\, x_{29}=7,\, x_{37}=-7,\, x_{41}=6,\, x_{53}=3,\, x_{61}=15.$$
\end{remark}

\begin{conjecture} Let $p>3$ be a prime.

{\rm (i)} If $p\eq1\pmod4$, then
\begin{equation}\l|\l(\f{j+k}p\r)-\l(\f{j-k}p\r)\r|_{1\ls j,k\ls (p-1)/2}
=(-p)^{(p-1)/4}
\end{equation}
and
\begin{equation}\begin{aligned}&\ \l|x+\l(\f{j+k}p\r)+\l(\f{j-k}p\r)+\l(\f jp\r)y+\l(\f kp\r)z+\l(\f{jk}p\r)w\r|_{1\ls j,k\ls (p-1)/2}
\\=&\ (-p)^{(p-5)/4}\l(\l(\f{p-1}2\r)^2wx-\l(\f{p-1}2y-1\r)\l(\f{p-1}2z-1\r)\r).
\end{aligned}
\end{equation}

{\rm (ii)} When $p\eq3\pmod4$, we have
\begin{equation}\begin{aligned}&\ (-1)^{\f{h(-p)+1}2}\l|x+\l(\f{j+k}p\r)+\l(\f{j-k}p\r)+\l(\f jp\r)y+\l(\f{jk}p\r)w\r|_{1\ls j,k\ls (p-1)/2}
\\=&\ p^{(p-3)/4}\l(\f{p-1}2y-1+\l(2-\l(\f 2p\r)\r)h(-p)(w+x)-\l(\f 2p\r)\f{16q_p}pwx\r)
\end{aligned}
\end{equation}
for some integer $q_p$ only depending on $p$.
\end{conjecture}
\begin{remark}
Our computation indicates that
\begin{gather*}q_7=q_{11}=1,\,q_{19}=9,\,q_{23}=15,\,q_{31}=24,\,q_{43}=27,\,q_{47}=72,\,q_{59}=62,
\\q_{67}=51,\,q_{71}=259,\,q_{79}=82,\,q_{83}=18,\,q_{103}=349,\,q_{107}=-68,\,q_{127}=478.
\end{gather*}
\end{remark}

\begin{conjecture} Let $p>3$ be a prime. If $p\eq1\pmod4$, then
\begin{equation}\begin{aligned}&\ \l|x+\l(\f{j+k}p\r)+\l(\f{j-k}p\r)+\l(\f jp\r)y+\l(\f kp\r)z+\l(\f{jk}p\r)w\r|_{0\ls j,k\ls (p-3)/2}
\\&\qquad=\l(\f 2p\r)p^{(p-5)/4}(px-wx+(y+2)(z+2)).
\end{aligned}
\end{equation}
When $p\eq3\pmod4$, there is an integer $m_p\in\Z$ only depending on $p$ such that
\begin{equation}\begin{aligned}&\ \l|x+\l(\f{j+k}p\r)+\l(\f{j-k}p\r)+\l(\f jp\r)y+\l(\f kp\r)z+\l(\f{jk}p\r)w\r|_{0\ls j,k\ls (p-3)/2}
\\&\qquad= (-1)^{(h(-p)+1)/2}p^{(p-7)/4}(p-2m_pw)x.
\end{aligned}
\end{equation}
\end{conjecture}
\begin{remark} Our computation indicates that
$$m_7=2,\, m_{11}=1,\, m_{19}=-3,\, m_{23}=-1,\, m_{31}=3,\, m_{43}=1,\, m_{47}=0,\, m_{59}=8.$$
\end{remark}

\section{Conjectures on determinants involving
  $(\f{j\pm k\pm1}p)$ or $(\f{j^2\pm k^2}p)$}
\setcounter{lemma}{0}
\setcounter{theorem}{0}
\setcounter{corollary}{0}
\setcounter{remark}{0}
\setcounter{conjecture}{0}
\setcounter{equation}{0}

\begin{conjecture} Let $p>3$ be a prime.

{\rm (i)} If $p\eq1\pmod4$, then
\begin{equation}\begin{aligned}&\ \l|x+\l(\f{j+k-1}p\r)+\l(\f jp\r)y+\l(\f kp\r)z+\l(\f{jk}p\r)w\r|_{1\ls j,k\ls(p-1)/2}
\\=&\ \l(\f 2p\r)2^{(p-1)/2}((yz-(w+1)x)a_p+(w(1-x)+(y+1)(z+1))b_p).
\end{aligned}
\end{equation}
When $p\eq3\pmod4$, we have
\begin{equation}\begin{aligned}&\ \l|x+\l(\f{j+k-1}p\r)+\l(\f jp\r)y+\l(\f kp\r)z+\l(\f{jk}p\r)w\r|_{1\ls j,k\ls(p-1)/2}
\\=&\ 2^{(p-1)/2}(yz-(w+1)x).
\end{aligned}\end{equation}

{\rm (ii)} If $p\eq1\pmod4$, then
\begin{equation}\begin{aligned}&\ \l|x+\l(\f{j+k-1}p\r)+\l(\f jp\r)y+\l(\f kp\r)z+\l(\f{jk}p\r)w\r|_{1\ls j,k\ls(p+1)/2}
\\=&\ \l(\f 2p\r)2^{(p-1)/2}(pb_p((w+1)x-yz)+a_p(w(x-1)-(y+1)(z+1)).
\end{aligned}\end{equation}
When $p\eq3\pmod4$, we have
\begin{equation}\begin{aligned}&\ \l|x+\l(\f{j+k-1}p\r)+\l(\f jp\r)y+\l(\f kp\r)z+\l(\f{jk}p\r)w\r|_{1\ls j,k\ls(p+1)/2}
\\=&\ 2^{(p-1)/2}(w(1-x)+(y+1)(z+1)).
\end{aligned}\end{equation}
\end{conjecture}
\begin{remark} In the case $w=0$, this reduces to Theorem \ref{Th1.4}.
\end{remark}

\begin{conjecture} \label{Conj-p-3} Let $p>3$ be a prime. If $p\eq1\pmod4$, then
\begin{equation}\begin{aligned}&\ \l|x+\l(\f{j+k-1}p\r)+\l(\f jp\r)y+\l(\f kp\r)z+\l(\f{jk}p\r)w\r|_{1\ls j,k\ls(p-3)/2}
\\=&\ \l(\f 2p\r)2^{(p-5)/2}\l(\l(2-\l(\f 2p\r)\r)a_p-pb_p\r)x
\\&\ +\l(\f 2p\r)2^{(p-5)/2}\l(a_p+\l(\l(\f 2p\r)-2\r)b_p\r)(w+y+z+1)
\\&\ +\l(\f 2p\r)2^{(p-5)/2}\l((p-1)b_p+\l(\l(\f 2p\r)-1\r)(a_p+b_p)\r)(yz-wx).
\end{aligned}\end{equation}
When $p\eq3\pmod4$, we have
\begin{equation}\begin{aligned}&\ \l|x+\l(\f{j+k-1}p\r)+\l(\f jp\r)y+\l(\f kp\r)z+\l(\f{jk}p\r)w\r|_{1\ls j,k\ls(p-3)/2}
\\=&\ 2^{(p-5)/2}\l(\l(\l(\f 2p\r)-2\r)x-w-y-z-1+\l(1-\l(\f 2p\r)\r)(yz-wx)\r).
\end{aligned}\end{equation}
\end{conjecture}
\begin{remark} Let $p=2n+1>3$ be a prime. As
\begin{align*}\l|x+\l(\f{j+k-1}p\r)\r|_{1\ls j,k\ls n-1}
&=\l|x+\l(\f{n+1-s+(n+1-t)-1}p\r)\r|_{2\ls s,t\ls n}
\\&= (-1)^{n(n-1)}\l|(-1)^nx+\l(\f{j+k}p\r)\r|_{2\ls j,k\ls n},
\end{align*}
Conjecture \ref{Conj-p-3} implies that
\begin{equation}\begin{aligned}&\ \l|x+\l(\f{j+k}p\r)\r|_{2\ls j,k\ls(p-1)/2}
\\=&\ \begin{cases}(\f 2p)2^{(p-5)/2}(a_p-pb_px+(2-(\f 2p))(a_px-b_p))&\t{if}\ p\eq1\pmod4,
\\2^{(p-5)/2}((2-(\f 2p))x-1)&\t{if}\ p\eq3\pmod4.\end{cases}\end{aligned}
\end{equation}
\end{remark}

\begin{conjecture}\label{Conj4.2} Let $p>3$ be a prime. If $p\eq1\pmod4$, then
\begin{equation}\begin{aligned}&\ \l|x+\l(\f{j+k-1}p\r)+\l(\f jp\r)y+\l(\f kp\r)z+\l(\f{jk}p\r)w\r|_{0\ls j,k\ls(p-1)/2}
\\=&\ \l(\f 2p\r)2^{(p-3)/2}(pb_p-2a_p)(w+y+z+1)
\\&\ +\l(\f 2p\r)2^{(p-3)/2}\l((2b_p-a_p)px+((p-2)a_p-pb_p)(yz-wx)\r).
\end{aligned}\end{equation}
If $p\eq3\pmod4$, then
\begin{equation}\begin{aligned}&\ \l|x+\l(\f{j+k-1}p\r)+\l(\f jp\r)y+\l(\f kp\r)z+\l(\f{jk}p\r)w\r|_{0\ls j,k\ls(p-1)/2}
\\=&\ 2^{(p-3)/2}\l(2w(1-x)+2(y+1)(z+1)+p((w+1)x-yz)\r).
\end{aligned}\end{equation}
\end{conjecture}
\begin{remark} Let $p=2n+1>3$ be a prime. As
\begin{align*}\l|x+\l(\f{j+k-1}p\r)\r|_{0\ls j,k\ls n}
=&\l|x+\l(\f{n+1-s+(n+1-t)-1}p\r)\r|_{1\ls s,t\ls n+1}
\\=&\ (-1)^{n(n+1)}\l|(-1)^nx+\l(\f{j+k}p\r)\r|_{1\ls j,k\ls n+1},
\end{align*}
Conjecture \ref{Conj4.2} implies that
\begin{equation}\begin{aligned}&\ \l|x+\l(\f{j+k}p\r)\r|_{1\ls j,k\ls(p+1)/2}
\\=&\ \begin{cases}(\f 2p)2^{(p-3)/2}(pb_p-2a_p+(2b_p-a_p)px)&\t{if}\ p\eq1\pmod4,
\\2^{(p-3)/2}(2-px)&\t{if}\ p\eq3\pmod4.\end{cases}\end{aligned}
\end{equation}
\end{remark}

\begin{conjecture} \label{Conj4.3} For any prime $p>3$ with $p\eq3\pmod4$, we have
\begin{equation}\begin{aligned}&\ \l|x+\l(\f{j+k-1}p\r)+\l(\f jp\r)y+\l(\f kp\r)z+\l(\f{jk}p\r)w\r|_{0\ls j,k\ls(p-3)/2}
\\=&\ 2^{(p-5)/2}\l(\l(\f p2\l(2-\l(\f 2p\r)\r)-4\r)x+\l(\f p2-2-\l(\f 2p\r)\r)(w+y+z+1)\r)
\\&\ +2^{(p-5)/2}\l(\f p2\l(\l(\f 2p\r)-1\r)+2-\l(\f 2p\r)\r)(yz-wx).
\end{aligned}
\end{equation}
\end{conjecture}
\begin{remark} Let $p=2n+1>3$ be a prime. As
\begin{align*} \l|x+\l(\f{j+k-1}p\r)\r|_{0\ls j,k\ls n-1}
=& \l|x+\l(\f{n+1-s+(n+1-t)-1}p\r)\r|_{2\ls s,t\ls n+1}
\\=& (-1)^n\l|(-1)^nx+\l(\f{j+k}p\r)\r|_{2\ls j,k\ls n+1},
\end{align*}
when $p\eq3\pmod4$ Conjecture \ref{Conj4.3} implies that
\begin{equation}\begin{aligned}&\ \l|x+\l(\f{j+k}p\r)\r|_{2\ls j,k\ls(p+1)/2}
\\=&\ 2^{(p-5)/2}\l(2+\l(\f 2p\r)-\f p2+\l(\f p2\l(2-\l(\f 2p\r)\r)-4\r)x\r).\end{aligned}
\end{equation}
\end{remark}

\begin{conjecture} Let $p>3$ be a prime. If
$p\eq1\pmod4$, then
\begin{equation}\begin{aligned}&\ \l|x+\l(\f{j+k+1}p\r)+\l(\f jp\r)y+\l(\f kp\r)z+\l(\f{jk}p\r)w\r|_{0\ls j,k\ls(p-1)/2}
\\=&\ \l(\f 2p\r)2^{(p-1)/2}pb_p\l(x+\f{p-2}2(yz-wx)\r)
\\&\ +\l(\f 2p\r)2^{(p-1)/2}a_p\l(w\l( x+\f{p-2}2\r)-(y+1)(z+1)\r).
\end{aligned}
\end{equation}
When $p\eq3\pmod4$, we have
\begin{equation}\begin{aligned}&\ \l|x+\l(\f{j+k+1}p\r)+\l(\f jp\r)y+\l(\f kp\r)z+\l(\f{jk}p\r)w\r|_{0\ls j,k\ls(p-1)/2}
\\=&\ 2^{(p-1)/2}\l(w\l(\f{p-2}2-x\r)+(y+1)(z+1)\r).
\end{aligned}
\end{equation}
\end{conjecture}
\begin{remark} For any prime $p=2n+1>3$, by \eqref{C*} we have
\begin{align*}&\ \l|x+\l(\f{j+k+1}p\r)\r|_{0\ls j,k\ls n}=\l|x+\l(\f{(j+1)+(k+1)-1}p\r)\r|_{0\ls j,k\ls n}
\\=&\ \det C_p^*(x)=\begin{cases}(-1)^{n/2}2^n(pb_px-a_p)&\t{if}\ p\eq1\pmod4,\\2^n&\t{if}\ p\eq3\pmod4.
\end{cases}
\end{align*}
\end{remark}

\begin{conjecture} Let $p>3$ be a prime. If
$p\eq1\pmod4$, then
\begin{equation}\begin{aligned}&\ \l|x+\l(\f{j+k+1}p\r)+\l(\f jp\r)y+\l(\f kp\r)z\r|_{1\ls j,k\ls(p-1)/2}
\\=&\ \l(\f 2p\r)2^{(p-3)/2}\l((pb_p-2a_p)x+2(n_p+b_p-a_p)yz\r)
\\&\ +\l(\f 2p\r)2^{(p-3)/2}\l((2b_p-a_p-1)(y+z+1)+1\r)
\end{aligned}
\end{equation}
for some positive integer $n_p$.
When $p\eq3\pmod4$, we have
\begin{equation}\begin{aligned}&\ \l|x+\l(\f{j+k+1}p\r)+\l(\f jp\r)y+\l(\f kp\r)z+\l(\f{jk}p\r)w\r|_{1\ls j,k\ls(p-1)/2}
\\=&\ 2^{(p-3)/2}\l(1-(-1)^{(h(-p)-1)/2}\r)(y+z+2(yz-wx))
\\&\ +2^{(p-3)/2}\l((p-3)\l(yz-wx+\f w2\r)-2x+1\r).
\end{aligned}
\end{equation}
\end{conjecture}
\begin{remark} Our computation indicates that
$$n_5=1,\ n_{13}=11,\ n_{17}=39,\ n_{29}=68,\ n_{37}=230,\ n_{41}=1441,\ n_{53}=256.$$
For any odd prime $p=2n+1$, we clearly have
\begin{align*}\l|x+\l(\f{j+k+1}p\r)\r|_{1\ls j,k\ls n}=&\l|x+\l(\f{(n-j)+(n-k)+1}p\r)\r|_{0\ls j,k\ls n-1}
\\=&\ (-1)^n\l|(-1)^nx+\l(\f{j+k}p\r)\r|_{0\ls j,k\ls n-1}.
\end{align*}
\end{remark}

\begin{conjecture}\label{Conj-j-k+10} Let $p>3$ be a prime. If $p\eq1\pmod4$, then
\begin{equation}\begin{aligned}&\ \l|x+\l(\f{j-k+1}p\r)+\l(\f jp\r)y+\l(\f kp\r)z+\l(\f{jk}p\r)w\r|_{0\ls j,k\ls(p-1)/2}
\\=&\ (pb_p'-a_p')(w(1-x)+(y+1)(z+1))
\\&\ +p\l(wx-(y+1)z+\l(\f 2p\r)(b_p'-a_p')((1+w)x-yz)\r).
\end{aligned}
\end{equation}
When $p\eq3\pmod 4$, we have
\begin{equation}\begin{aligned}&\ \l|x+\l(\f{j-k+1}p\r)+\l(\f jp\r)y+\l(\f kp\r)z+\l(\f{jk}p\r)w\r|_{0\ls j,k\ls(p-1)/2}
\\=&\ 1-\l(\f 2p\r)px+w+y+\l(p\l(\f 2p\r)(-1)^{(h(-p)-1)/2}-1\r)z
\\&\ +\l(p\l(\f 2p\r)\l(1+(-1)^{(h(-p)-1)/2}\r)-1\r)(yz-wx).
\end{aligned}
\end{equation}
\end{conjecture}

\begin{conjecture}\label{Conj-j-k+1} Let $p>3$ be a prime.
If $p\eq1\pmod4$, then
\begin{equation}\begin{aligned}&\ \l|x+\l(\f{j-k+1}p\r)+\l(\f jp\r)y+\l(\f kp\r)z+\l(\f{jk}p\r)w\r|_{0\ls j,k\ls(p-3)/2}
\\=&\ (pb_p'-a_p')((w+1)x-yz)+\l(\f 2p\r)(wx-(y+1)z)
\\&\ +\l(\f 2p\r)(b_p'-a_p')(w(1-x)+(y+1)(z+1)).
\end{aligned}
\end{equation}
When $p\eq3\pmod4$, we have
\begin{equation}\begin{aligned}&\ \l|x+\l(\f{j-k+1}p\r)+\l(\f jp\r)y+\l(\f kp\r)z+\l(\f{jk}p\r)w\r|_{0\ls j,k\ls(p-3)/2}
\\=&\ x-\l(\f 2p\r)(w+y-z+1)-(-1)^{(h(-p)-1)/2}z
\\&\ +\l(1+(-1)^{(h(-p)-1)/2}-\l(\f 2p\r)\r)(wx-yz).
\end{aligned}
\end{equation}
\end{conjecture}
\begin{conjecture} Let $p>3$ be a prime. If $p\eq1\pmod4$, then
\begin{equation}\begin{aligned}&\ \l|x+\l(\f{j-k+1}p\r)+\l(\f jp\r)y+\l(\f kp\r)z+\l(\f{jk}p\r)w\r|_{1\ls j,k\ls(p-1)/2}
\\=&\ \l(\f 2p\r)\f{p-1}2((y+1)z-wx)+(pb_p'-a_p')((w+1)x-yz)
\\&\ +\l(\f 2p\r)(b_p'-a_p')(w(1-x)+(y+1)(z+1)).
\end{aligned}
\end{equation}
When $p\eq3\pmod4$, we have
\begin{equation}\begin{aligned}&\ \l|x+\l(\f{j-k+1}p\r)+\l(\f jp\r)y+\l(\f kp\r)z+\l(\f{jk}p\r)w\r|_{1\ls j,k\ls(p-1)/2}
\\=&\ (w+1)x-yz+\l(\f 2p\r)((y+1)(z-1)-w(x+1))
\\&\ +(-1)^{(h(-p)-1)/2}\,\f{p+1}2(wx-(y+1)z).
\end{aligned}
\end{equation}
\end{conjecture}

\begin{conjecture} Let $p$ be an odd prime.

{\rm (i)} When $p\eq1\pmod4$, for any $\da_1,\da_2\in\{\pm1\}$ the number
$$2\l|\l(\f{j+k}p\r)+\l(\f{j-k}p\r)+\da_1\l(\f{j^2+\da_2 k^2}p\r)\r|_{0\ls j,k\ls(p-1)/2}$$
is a quadratic residue modulo $p$.

{\rm (ii)} If $p\eq3\pmod4$, then the number
$$2\l|\l(\f{j+k}p\r)+\l(\f{j^2+k^2}p\r)\r|_{0\ls j,k\ls(p-1)/2}$$
is a quadratic residue modulo $p$.
\end{conjecture}

\section{Conjectures on determinants of the form
  $\{c,d\}_n=|(\f{j^2+cjk+dk^2}n)|_{1<j,k<n-1}$}
\setcounter{lemma}{0}
\setcounter{theorem}{0}
\setcounter{corollary}{0}
\setcounter{remark}{0}
\setcounter{conjecture}{0}
\setcounter{equation}{0}

Let $p$ be a prime with $p\eq1\pmod4$. By a classical result conjectured by Fermat and confirmed by Euler, $p=a^2+b^2$ for some $a,b\in\Z$ with $a$ odd. In view of Jacobsthal's theorem (cf. Theorem 6.2.9 of \cite[p.\,195]{BEW}),
$$p=\(\sum_{x=1}^{(p-1)/2}\l(\f{x(x^2+1)}p\r)\)^2+\(\sum_{x=1}^{(p-1)/2}\l(\f{x(x^2+d)}p\r)\)^2$$
for any $d\in\Z$ with $(\f dp)=-1$. As $x^2\eq-1\pmod p$ for a unique number $x\in\{1,\ldots,(p-1)/2\}$, we have
$$\f12\sum_{x=1}^{p-1}\l(\f{x(x^2+1)}p\r)=\sum_{x=1}^{(p-1)/2}\l(\f{x(x^2+1)}p\r)=\pm a\not=0.$$
Motivated by this, we obtain the following result.

\begin{theorem} Let $n>1$ be an integer with $n\eq1\pmod4$. Then $n$ is not a sum of two squares
if and only if
$$\sum_{x=0}^{n-1}\l(\f {x(x^2+1)}n\r)=0.$$
\end{theorem}
\Proof. Write $n=p_1^{a_1}\cdots p_k^{a_k}$, where $p_1,\ldots,p_k$ are distinct odd primes and $a_1,\ldots,a_k\in\N$. Applying the Chinese Remainder Theorem, we obtain
\begin{align*} \sum_{x=0}^{n-1}\l(\f{x(x^2+1)}n\r)
&=\sum_{x_1=0}^{p_1^{a_1}-1}\cdots\sum_{x_k=0}^{p_k^{a_k}-1}\prod_{s=1}^k\l(\f{x_s(x_s^2+1)}{p_s^{a_s}}\r)
\\&=\prod_{s=1}^k\sum_{x_s=0}^{p_s^{a_s}-1}\l(\f{x_s(x_s^2+1)}{p_s}\r)^{a_s}
\\&=\prod_{s=1}^k\sum_{q_s=0}^{p_s^{a_s-1}-1}\sum_{r_s=0}^{p_s-1}\l(\f{(p_sq_s+r_s)((p_sq_s+r_s)^2+1)}{p_s}\r)^{a_s}
\\&=\prod_{s=1}^k p_s^{a_s-1}\sum_{x=1}^{p_s-1}\l(\f{x(x^2+1)}{p_s}\r)^{a_s}.
\end{align*}
If $p_s\eq1\pmod4$, then $|\{1\ls x\ls p-1:\ x^2\eq-1\pmod{p_s}\}|=2$ and hence
$$\sum_{x=1}^{p_s-1}\l(\f{x(x^2+1)}{p_s}\r)^{a_s}\eq p_s-3\not\eq0\pmod {4}.$$
When $p_s\eq3\pmod4$ and $2\mid a_s$, we have
$$\sum_{x=1}^{p_s-1}\l(\f{x(x^2+1}p\r)^{a_s}=p_s-1\not=0.$$
When $p_s\eq3\pmod4$ and $2\nmid a_s$, we have
$$\sum_{x=1}^{p_s-1}\l(\f{x(x^2+1}p\r)^{a_s}=\sum_{x=1}^{(p_s-1)/2}\(\l(\f{x(x^2+1)}{p_s}\r)
+\l(\f{-x(x^2+1)}{p_s}\r)\)=0.$$

In view of the above, we see that $\sum_{x=0}^{n-1}(\f {x(x^2+1)}n)=0$ if and only if
$p_s\eq3\pmod4$ and $2\nmid a_s$ for some $s=1,\ldots,k$.
By a known result (cf. \cite[p.\,279]{IR}), $n$ is not a sum of two squares if and only if for some prime divisor $p$ of $n$
we have $p\eq3\pmod4$ and $\ord_p(n)\eq1\pmod2$, where $\ord_p(n)$ is the $p$-adic order of $n$ at the prime $p$.
So the desired result follows. \qed

Let $n>1$ be an odd integer, and let $c,d\in\Z$. The author \cite{S19} investigated the new kinds of determinants
$$[c,d]_n=\l|\l(\f{j^2+cjk+dk^2}n\r)\r|_{0\ls j,k\ls n-1}$$
and
$$(c,d)_n=\l|\l(\f{j^2+cjk+dk^2}n\r)\r|_{1\ls j,k\ls n-1},$$
where $(\f{\cdot}n)$ denotes the Jacobi symbol. Some conjectures on such determinants
were later confirmed by D. Krachun, F. Petrov, Z.-W. Sun and M. Vsemirnov \cite{KPSV}.
Now we introduce the new determinant
\begin{equation}\{c,d\}_n=\l|\l(\f{j^2+cjk+dk^2}n\r)\r|_{1<j,k<n-1}.
\end{equation}

\begin{conjecture} {\rm (i)} For any positive integer $n\eq1\pmod4$ which is not a sum of two squares,
we have $\{3,2\}_n=0$.

{\rm (ii)} For any positive integer $n\eq3\pmod 4$, we have $\f{\varphi(n)}2\mid \{3,2\}_n$, where $\varphi$ is Euler's totient function.

{\rm (iii)} For any positive integer $n\eq3\pmod8$, we have
$$\{3,2\}_n=\f{\varphi(n)}2x^2$$
for some $x\in\Z$.
\end{conjecture}
\begin{remark} We have verified this for all positive odd integers $n<2000$.
By \cite[Corollary 1.1]{KPSV}, $(3,2)_n=[3,2]_n=0$ for any positive integer $n\eq3\pmod4$.
\end{remark}

\begin{conjecture} {\rm (i)} We have $\{2,2\}_p=0$ for any prime $p\eq13,19\pmod{24}$.

{\rm (ii)} We have $\{2,2\}_p\eq0\pmod p$ for any prime $p\eq17,23\pmod{24}$.
\end{conjecture}
\begin{remark} We have verified this conjecture for odd primes $p<2000$.
\end{remark}

\begin{conjecture} {\rm (i)} We have $\{4,2\}_n=\{8,8\}_n=0$ for any positive integer $n\eq5\pmod8$.

{\rm (ii)} We have $\{3,3\}_n=0$ for any positive integer $n\eq5\pmod{12}$.
\end{conjecture}
\begin{remark} We have verified this conjecture for positive odd integers $n<2000$.
By \cite[Corollary 1.1]{KPSV}, $(4,2)_n=(8,8)_n=(3,3)_n=0$ for any positive integer $n\eq3\pmod4$.
\end{remark}

\begin{conjecture} We have $\{42,-7\}_n=\{21,112\}_n=0$ for any positive integer $n\eq1\pmod4$ with $(\f n7)=-1$.
\end{conjecture}
\begin{remark} We have verified this conjecture for positive odd integers $n<2000$.
By \cite[Theorem 1.1(iv)]{KPSV}, $(42,-7)_n=(21,112)_n=0$ for any positive integer $n$ with $(\f n7)=-1$.
\end{remark}

\begin{conjecture} {\rm (i)} Let $n>3$ be an odd integer.
Then $\{2,3\}_n\eq0\pmod n$. Moreover, $\{2,3\}_n\eq0\pmod{n^2}$ if $n\not\eq\pm1\pmod {12}$.

{\rm (ii)} For any odd integer $n>7$, we have $\{6,15\}_n\eq0\pmod{n}.$
\end{conjecture}
\begin{remark} We have verified this conjecture for positive odd integers $n<2000$.
The author \cite[Conjecture 4.8]{S19} conjectured that $(2,3)_n\eq0\pmod{n^2}$
for each odd integer $n>3$, and that $(6,15)_n\eq0\pmod{n^2}$ for any odd integer $n>5$.
\end{remark}

\begin{conjecture} {\rm (i)} For any positive integer $n\eq13,17\pmod{20}$
which is a sum of two squares, we have $\{5,5\}_n=0$.

{\rm (ii)} Let $n>1$ be an odd integer. We have $(\f{\{5,5\}_n}n)=0$ if $n\eq11,19\pmod{20}$,
or $n\eq9\pmod{60}$ and $n>69$.
\end{conjecture}
\begin{remark} We have verified this conjecture for positive odd integers $n<2000$.
By \cite[Theorem 1.4]{KPSV}, $[5,5]_p=0$ for any prime $p\eq13,17\pmod{20}$.
\end{remark}

\begin{conjecture} {\rm (i)} For any positive integer $n\eq5\pmod{12}$
which is a sum of two squares, we have $\{10,9\}_n=0$.

{\rm (ii)} We have $\{10,9\}_p\eq0\pmod p$ for any prime $p\eq11\pmod{12}$.
\end{conjecture}
\begin{remark} We have verified this conjecture for $n,p<2000$. By \cite[Theorem 1.4]{KPSV}, $(10,9)_p=0$ for any prime $p\eq5\pmod{12}$.
\end{remark}

\begin{conjecture} {\rm (i)} For any positive integer $n\eq13,17\pmod{24}$
which is a sum of two squares, we have $\{8,18\}_n=0$.

{\rm (ii)} We have $\{8,18\}_p\eq0\pmod {p^2}$ for any prime $p\eq19\pmod{24}$.

{\rm (iii)} We have $\{8,18\}_p\eq0\pmod {p}$ for any prime $p\eq23\pmod{24}$.
\end{conjecture}
\begin{remark} We have verified this conjecture for $n,p<2000$. In 2018, the author \cite{MF} conjectured
that $[8,18]_p=0$ for any prime $p\eq13,17\pmod{24}$, which was confirmed by Michael Stoll (cf. \cite{MF}) by using advanced tools such as elliptic curves with complex multiplication by $\Z[\sqrt{-6}]$ and $\ell$-adic Tate modules.
\end{remark}

\setcounter{conjecture}{0} 
\begin{thebibliography}{S19}

\bibitem{BEW} B. C. Berndt, R. J. Evans and K. S. Williams,
Gauss and Jacobi Sums, John Wiley \& Sons, 1998.


\bibitem {Ch} R. Chapman, {\it Determinants of Legendre symbol matrices}, Acta Arith.
{\bf 115} (2004), 231--244.

\bibitem{IR} K. Ireland and M. Rosen, {\it A Classical Introduction to Modern Number Theory},
2nd Edition, Grad. Texts. Math., vol. 84, Springer, New York, 1990.


\bibitem{KPSV} D. Krachun, F. Petrov, Z.-W. Sun and M. Vsemirnov, {\it On some determinants involving Jacobi symbols}, Finite Fields Appl. {\bf 64} (2020), Article 101672.
    
\bibitem{K05} C. Krattenthaler, {\it Advanced determinant calculus: a complement}, Linear Algebra Appl. {\bf 411} (2005), 68--116.

\bibitem{L} D. H. Lehmer, {\it On certain character matrices}, Pacific J. Math. {\bf 6} (1956), 491--499.

\bibitem{MF} Z.-W. Sun, {\it A series of conjectures on $\sum_{x=0}^{(p-1)/2}(\f{x^5+cx^3+dx}p)$ (III)},
Question 319259 at MathOverflow (with an answer from Michael Stoll), Dec. 22, 2018. https://mathoverflow.net/questions/319259

\bibitem{S19} Z.-W. Sun, {\it On some determinants with Legendre symbol entries}, Finite Fields Appl. {\bf 56} (2019), 285--307.

\bibitem{S19p} Z.-W. Sun, {\it Quadratic residues and related permutations and identities},
Finite Fields Appl. {\bf 59} (2019), 246--283.

\bibitem{RJ} Z.-W. Sun, {\it On some determinants involving the tangent functions}, Ramanujan J.
{\bf 64} (2024), 309--332.


\bibitem{V12} M. Vsemirnov, {\it On the evaluation of R. Chapman's ``evil determinant"},
Linear Algebra Appl. {\bf 436} (2012), 4101--4106.

\bibitem{V13} M. Vsemirnov, {\it On R. Chapman's ``evil determinant": case $p\eq1\ (\mo\ 4)$},
Acta Arith. {\bf 159} (2013),  331--344.

\bibitem{WWN} L.-Y. Wang, H.-L. Wu and H.-X. Ni,
{\it On a generalization of R. Chapman's ``evil determinant"}, preprint, arXiv:2405.02112, 2024.


\end{thebibliography}
\end{document}